\documentclass[12pt]{article}

\usepackage[latin1]{inputenc}
\usepackage{graphicx}
\usepackage{epsfig}

\newtheorem{theorem}{Theorem}
\newtheorem{lemma}{Lemma}
\newtheorem{remark}{Remark}

\newfont{\bb}{msbm10 at 12pt}
\def\r{\hbox{\bb R}}

\newenvironment{proof}{\trivlist
\item[\hskip\labelsep{\it Proof}\,:]}{\hfill{$q.e.d$}\endtrivlist}

\title{A new proof of a characterization \\
of small spherical caps}
\author{Rafael L\'opez\thanks{Partially
supported by MEC-FEDER
 grant no. MTM2007-61775.}\\
Departamento de Geometr\'{\i}a y Topolog\'{\i}a\\
Universidad de Granada\\
18071 Granada\\
Spain\\
email: rcamino@ugr.es}

\date{}

\begin{document}

\maketitle


\begin{abstract} It is known  that planar disks and small spherical caps are the only
constant mean curvature graphs whose boundary is a round circle.
Usually, the proof invokes the Maximum Principle for elliptic
equations. This paper presents a new proof of this result motivated
by an article due to Reilly. Our proof
utilizes a flux formula for surfaces with constant mean curvature
together with integral equalities on the surface.
\end{abstract}

\section{Introduction and the result}

A surface in Euclidean space $\r^3$ with the property that its mean curvature is constant at each point is called a \emph{constant mean
curvature surface} or CMC surface for short.  Round spheres are closed CMC surfaces. Here by closed surface we mean compact and
without boundary surface. A famous
theorem due to Hopf  asserts that any closed CMC surface of genus
$0$ must be a round sphere \cite{ho}. Later,   Alexandrov proved in 1956 that any embedded closed
CMC surface in $\r^3$ must be a round sphere \cite{al}. For a long
time it was an open question whether or not spheres were the only
closed CMC surfaces in $\r^3$. If a such surface were to exist, it
would necessarily be a surface with self-intersections and higher
genus. In 1986, Wente found an immersed torus with constant mean
curvature \cite{we}. This discovery inspired a great deal of work in
the search for new examples of closed CMC surfaces. For readers
interested in the subject, we refer the recent survey \cite{ke} and
references therein.

 We now consider compact CMC surfaces with non-empty boundary. The
 simplest case for the boundary is a round circle. If $C$ is a
 circle of radius $r>0$, we consider $C$ in a sphere $S(R)$ of radius $R$,
 $R\geq r$.  The mean curvature
of $S(R)$ is $H=1/R$ with respect to the inward orientation. Then $C$ splits
$S(R)$ in two spherical caps with the
 same boundary $C$ and constant mean curvature $H$.  If $R=r$, both
 caps are hemispheres whereas  if $R>r$, there are two geometrically distinct caps which we call the
small and the big spherical cap. On the other hand,  the planar disk
bounded by $C$
 is a compact surface with constant mean curvature $H=0$. These
 surfaces are the only totally  umbilic compact CMC surfaces
 bounded by $C$.

In 1991,    Kapouleas
 found  other examples of CMC surfaces bounded by a circle
\cite{ka}. The surfaces that he obtained have higher genus and  self-intersections. Thus, one asks under what conditions a compact CMC  surface bounded by a circle is spherical.
Taking into account the theorems  of Hopf and Alexandrov for closed
 surfaces above cited, the natural  hypotheses to consider for surfaces bounded by a circle is that either  $S$ has the simplest  possible topology, that is, the topology of a disk, or that $S$ is  embedded. Surprisingly, we have

\begin{enumerate}
\item[] Conjeture 1. Planar disks and spherical caps are the only
compact CMC surfaces bounded by a circle that are \emph{topological
disks}.
\item[] Conjeture 2. Planar disks and spherical caps are the only
compact CMC surfaces bounded by a circle that are \emph{embedded}.
\end{enumerate}

This means that  our knowledge about the structure of the space of
CMC surfaces bounded by a circle is  quite limited and only several partial results  have been  obtained by different authors (we refer to \cite{ke} again). Of course,  the methods of proof
for the  Hopf and Alexandrov Theorems can not be applied  with complete success in the
context of a non-empty boundary. This fact, together the lack of
examples, suggests that although the problems in the non-empty boundary
case have the same flavor as in the closed one, the proofs are more
difficult.

 A partial answer  to the conjecture 2 is the following
\begin{theorem}[Alexandrov]\label{lopez:embedded} Let $C$ be a round circle in a plane $\Pi$
and let  $S$ be an embedded  compact CMC surface bounded by $C$. If
$S$ lies in one side of $\Pi$, then $S$ is a planar disk or a
spherical cap.
\end{theorem}
The extra hypothesis that we add is that $S$ \emph{lies on one side of
the plane containing the boundary}. Although Alexandrov did not
state this result, the  proof
 is accomplished  using the same technique that he used  in proving his theorem that was stated above:
 the so-called Alexandrov reflection method. Behind this method
  lies the classical  Maximum
 Principle  for  elliptic partial differential equations, together with the moving plane technique.
A particular case of this Theorem is the following. Given $H$ and a circle $C$, among
the two spherical caps bounded by $C$  with mean curvature $H$, only
the small one is a graph. Using this method, we characterize the small
spherical caps as

\begin{theorem}\label{lopez:t2} Let $C$ be a round circle in
a plane $\Pi$ and let  $S$ be a  compact CMC surface bounded by $C$.
If $S$ is a graph over $\Pi$, then $S$ is a planar disk or a small
spherical cap.
\end{theorem}

The purpose of this article is to give a new proof of this result and  that \emph{does not
involve the Maximum Principle}. This different approach, that is, avoiding the Maximum Principle, appeared
 in  the closed  case which motivated the
present work. In 1978, Reilly
 obtained  another proof of the Alexandrov
theorem for CMC closed surfaces without the use of the Maximum Principle  thanks to a combination of the Minkowski formulae with some new elegant arguments \cite{re}. In this sense,  a new elementary proof of Alexandrov's Theorem  due to Ros appears in \cite{ro}.

In the same spirit, we will use integral formulae together a
type of "flux formula". Moreover we will see in the next section how
the equation that characterizes a CMC surface can be expressed in terms of the
Laplace operator. This was already noticed by Reilly as one can see from  the title of
his article. This fact allows us to establish some geometrical
properties about CMC surfaces using basic multivariable Calculus.

\section{Mean curvature, graphs and the  Laplacian}\label{lopez:alexa}

Let $S$ be a surface in $\r^3$ and which we  write  locally as the
graph of a smooth function $f$, $z=f(x,y)$,  $(x,y,z)$ being the
usual coordinates of $\r^3$. We orient   $S$ with the choice of normal
given by
\begin{equation}\label{lopez:gauss}
N(x,y,f(x,y))= \frac{(-f_x,-f_y,1)}{\sqrt{1+f_x^2+f_y^2}}(x,y),
\end{equation} where the subscripts  indicate the corresponding
partial  derivatives.  The mean curvature $H$ of $S$  satisfies the
following partial differential equation:
\begin{equation}\label{lopez:media}
2H(1+f_x^2+f_y^2)^{ \frac{3}{2}}=(1+f_y^2)f_{xx}-2f_x f_y
f_{xy}+(1+f_x^2)f_{yy}.
\end{equation}
Equation (\ref{lopez:media}) may written  as
\begin{equation}\label{lopez:media2}
\mbox{div }\left(\frac{\nabla f}{\sqrt{1+|\nabla f|^2}}\right)=2H,
\end{equation}
where  div and $\nabla$ stand for the divergence and gradient operators
respectively. In  PDE theory, this equation falls into the category  of elliptic  type,
whose main property is the existence of a Maximum Principle: if two functions $f_1$
and $f_2$ satisfy equation (\ref{lopez:media2}) with the same
Dirichlet condition, then $f_1=f_2$. We refer the reader to \cite[sect. 10.5]{gt}.
For CMC surfaces, this
 geometrically translates to the assertion
that if two CMC surfaces with the same  constant mean curvature
intersect tangentially at some point and one surfaces lies locally on
one side of the other, then both surfaces must coincide in a
neighborhood of that point. This property was used by Alexandrov in proving
his theorem by using a surface as a comparison surface with itself in a
 reflection process. When the surface has non-empty boundary, we
must add the extra hypothesis that $S$ lies over $\Pi$ as
states Theorem \ref{lopez:embedded}. By doing this, we avoid the
presence of a possible contact  between an interior point with a
boundary point of the surface where the Maximum Principle fails. We
refer the reader to \cite{ke,ko} for detailed proofs of Theorem
\ref{lopez:embedded}.

We prove two results about CMC surfaces which have their own geometric  interest (they will not be used
later). Both results are well known in the literature although usually the Maximum Principle is invoked in the proofs. However, we show them  by using a basic knowledge of calculus and differential geometry.

\begin{theorem}\label{lopez:lado} Let $D$ be a domain of a plane $\Pi$ and
let $S$ be  a compact CMC graph on  $D$ whose boundary is $\partial
D$. If $H\not=0$, then $int(S)$ lies in one side of $\Pi$.
\end{theorem}
\begin{proof}
We  argue by contradiction. Assume that $S$ has (interior) points on
both sides of $\Pi$. Consider $p_0=(x_0,y_0,z_0)$, $z_0>0$ and
$p_1=(x_1,y_1,z_1)$, $z_1\leq 0$, points of $S$ with highest and
lowest height $z$, respectively. If $S$ is the graph of a  function
$z=f(x,y)$, then,
\begin{eqnarray}
& &f_x(x_i,y_i)=f_y(x_i,y_i)=0,\hspace*{.5cm}i=0,1\label{lopez:ayu1}\\
& & f_{xx}(x_0,y_0),f_{yy}(x_0,y_0)\leq 0,\hspace*{.5cm} f_{xx}(x_1
,y_1),f_{yy}(x_1 ,y_1)\geq 0.\label{lopez:ayu2}
\end{eqnarray}
Using  the
orientation given by (\ref{lopez:gauss}), let us  compute the mean
curvature $H$ at $p_0$ and $p_1$. Because  $H$ is constant and using
(\ref{lopez:ayu1}) and (\ref{lopez:ayu2}), equation
(\ref{lopez:media}) leads to
\begin{equation}\label{lopez:hache}
2H=2H(p_0)=\left( f_{xx}+f_{yy}\right)(x_0,y_0)\leq 0\leq \left(
f_{xx}+f_{yy}\right)(x_1,y_1)=2H(p_1)=2H. \end{equation}
 Since
$H\not=0$, we get a contradiction.
\end{proof}
The inequalities in (\ref{lopez:hache}) can written as $$2H(p_0)=\Delta_0
f(x_0,y_0)\leq 0\leq \Delta_0 f(x_1,y_1)=2H(p_1),$$
 where $\Delta_0=\partial_{xx}+\partial_{yy}$
is the Euclidean Laplacian. This indicates that under a certain choice
of coordinates,  (\ref{lopez:media})  can be expressed in terms of
the Laplace operator. See also \cite{re}.

 We treat the minimal case, that is, $H=0$.

\begin{theorem}
 Consider a Jordan curve  $C$  in a plane $\Pi$.
If $S$ is  a compact CMC surface with  $H\equiv 0$, whose boundary is $C$, then
$S$ is the planar domain $D$ that bounds $C$.
\end{theorem}

\begin{proof}  We point out that we have dropped the
hypothesis that $S$ is a graph. We use a similar proof as in Theorem
\ref{lopez:lado} and we follow the notation used there. The
reasoning is by contradiction again. Without loss of generality, we
assume  that $S$ has points over $\Pi$. Let $\Gamma\subset\Pi$ be a
circle of radius $r$ sufficiently large so that $D$ lies strictly
inside of the circular disk determined by $\Gamma$ and so that the
hemisphere $K$ with $\partial K=\Gamma$  over $\Pi$ also lies
over $S$. Let $S(H)$ be the family of small spherical caps  over
$\Pi$ with $\partial S(H)=\Gamma$ and parameterized by their mean
curvature $H$  oriented by (\ref{lopez:gauss}). Then
$-1/r<H<0$. In the limit case, $S(-1/r)=K$. Beginning from the value
$H=-1/r$, we let $H\rightarrow 0$ until the first value of $h$,
$-1/r<h<0$, that $S(h)$ touches the original surface $S$. Let $p_0$
be the contact point.  Both surfaces are
locally graphs of functions defined in the (common) tangent plane at
$p_0$. This point is not necessarily  the highest point of $S$,
but we do a change of coordinates so that $p_0$ is the origin, the
tangent plane of $S(h)$ and $S$ at $p_0$ is the $xy$-plane and
$S(h)$ lies over $S$ in a neighborhood of $p_0$. Now $p_0$ is the
highest point of $S$ and both surfaces lie below $\Pi$.

 Consider  the two functions $f$ and $g$ whose graphs
are $S$ and $S(h)$ respectively and defined in some planar domain
$\Omega$ of $\Pi$ containing the origin.  Let $u=f-g$. Then $u\leq 0$
on $\Omega$ with a local maximum at $(0,0)$. Because $g>0$ on
$\partial D$, $p_0$ is an interior point of $\Omega$. Consequently,
$f_x(0,0)=f_y(0,0)=g_x(0,0)=g_y(0,0)=0$ and
\begin{equation}\label{lopez:uxx}
u_{xx}(0,0),\ \ u_{yy}(0,0)\leq 0. \end{equation}
 However at the point $p_0$, equation
(\ref{lopez:media}) for $S$ and $S(h)$ is
$$0=(f_{xx}+f_{yy})(0,0)\hspace*{.5cm}\mbox{and}\hspace*{.5cm}
-\frac{1}{r}=(g_{xx}+g_{yy})(0,0),$$
respectively. By substracting
both equations, we obtain $1/r=(u_{xx}+u_{yy})(0,0)>0$,
contradicting (\ref{lopez:uxx}).

\end{proof}

\section{The effect of the boundary in the shape of a CMC surface}\label{lopez:effect}

We have seen that if $C$ is a circle of radius $r$, the possible
values of the mean curvatures $H$ for spherical caps bounded by $C$ lies
in the range $[-1/r,1/r]$ because $R\geq r$ for the radius of the
spheres $S(R)$. Thus, the boundary $C$ imposes restrictions on the
possible values of mean curvature. We show that this occurs for a
general curved boundary.   Consider  a compact CMC surface $S$ with
boundary $\partial S=C$ and let $Y$ be a variation field in
$\r^3$. The first variation formula of the area $|A|$ of the surface
$S$ along $Y$ is
$$\delta_Y|A|=-2 H\int_S \langle N,Y\rangle\ {\rm d}S-
\int_{\partial S}\langle \nu,Y\rangle\ {\rm d}s,$$ where $N$ is the unit normal vector of $S$, $H$ is the mean curvature relative to $N$,
$\nu$  represents  the inward unit vector along $\partial S$ and
${\rm d}s$ is the arc-length element of $\partial S$. Let us fix a
vector $\vec{a}\in\r^3$ and consider $Y$ as the generating field of a family of
translations in the direction of $\vec{a}$. As $Y$ generates
isometries of $\r^3$, the first variation of $A$ vanishes and thus
\begin{equation}\label{lopez:ba}
2H\int_S \langle N,\vec{a}\rangle\ {\rm d}S+\int_{\partial S}\langle
\nu,\vec{a}\rangle\ {\rm d}s=0.
\end{equation}
The first integral transforms into an integral over the boundary as
follows. The divergence of the field $ Z_p= (p\times \vec{a})\times
N$, $p\in S$, is $-2\langle N,\vec{a}\rangle$ (here $\times$ denotes
the cross product of $\r^3$). The Divergence Theorem, together with
(\ref{lopez:ba}), yields
\begin{equation}\label{lopez:flux}
-H\int_{\partial S}\langle\alpha\times \alpha',\vec{a}\rangle \ {\rm
d}s=\int_{\partial S}\langle \nu,\vec{a}\rangle\ {\rm d}s,
\end{equation}
where $\alpha$ is a parametrization of $\partial S$ such that
$\alpha'\times\nu=N$.

This equation  known as the ``balancing formula" or ``flux formula" is  due to R. Kusner
(see \cite{ku1}; also in \cite{ku2,kks}). It is a conservation law in the sense of
Noether that reflects the fact that the area (the potential) is
invariant under the group of translations of Euclidean space. On the
other hand,  if $D$ is a $2$-cycle with boundary $\partial S$,
the formula can be viewed as expressing the physical equilibrium between the force of
exterior pressure acting  on  $D$ (the left-hand side in
(\ref{lopez:flux})) with the force of  surface tension of $S$
that act along its boundary (the right-hand side).

If the boundary $C$  lies in the plane $\Pi=\{x\in\r^3;\langle
x,\vec{a}\rangle=0\}$, for $|\vec{a}|=1$, then (\ref{lopez:flux}) gives
\begin{equation}\label{lopez:ba2}
2H \bar{A}=\int_{\partial S}\langle\nu,\vec{a}\rangle\ {\rm d}s,
\end{equation}
where $\bar{A}$ is the algebraic area of $C$.   Given a closed curve
$C\subset\r^3$ that bounds a domain $D$, and noting that
$\langle\nu,\vec{a}\rangle\leq 1$, the  possible values of the mean curvature $H$  of $S$  satisfy
\begin{equation}\label{lopez:nece}
|H|\leq\frac{{\rm length}(C)}{2\ {\rm area}(D)}.
\end{equation}
 In  particular, {\it if $C$ is a
circle of radius $r>0$, a necessary condition for the existence of a
surface spanning $C$  with constant mean curvature $H$ is that
$|H|\leq 1/r$}.

\begin{remark}\label{lopez:gra} If $S$ is the graph of $z=f(x,y)$
it follows from the Divergence Theorem and (\ref{lopez:media2}) that
\begin{eqnarray*}2|H|{\rm area}(D)&=&\left| \int_D 2H\ {\rm d} D
\right|=\left|\int_{\partial D}
\langle \frac{\nabla f}{\sqrt{1+|\nabla f|^2}},\vec{n}\rangle\ {\rm d}s\right|\\
&\leq&\int_{\partial D}\frac{|\nabla f|}{\sqrt{1+|\nabla f|^2}}\
{\rm d}s<\int_{\partial D} 1\ {\rm d}s= {\rm
length}(C),\end{eqnarray*} where $\vec{n}$ is the unit normal vector
to  $\partial D$ in $\Pi$. Then,
$$|H|<\frac{\rm{length}(C)}{2\ \rm{ area}(D)}.$$
\end{remark}

As consequence of  Remark \ref{lopez:gra}, the proof of Theorem \ref{lopez:t2}
 is  very simple by using the Maximum Principle as we  show at this time.
 If $S$ is
the graph of a function $z=f_1(x,y)$ and  $H$ is its mean curvature,
then $|H|<1/r$, where $r$ is the radius of $C$. But there exists a
small spherical cap with the same boundary and mean curvature as
$S$. As this cap is the graph of  a function $f_2$, we have  that
$f_1$ and $f_2$ are two solutions of (\ref{lopez:media2}) with the
same Dirichlet condition, the Maximum Principle implies $f_1=f_2$. Other proof of Theorem \ref{lopez:t2} using a combination of the flux formula and the Maximum Principle appears in \cite{ba}.

\section{A new proof of Theorem \ref{lopez:t2}}
In this section we will prove our result without the use of the
 Maximum Principle. Let $S$  satisfy  the hypotheses of  Theorem \ref{lopez:t2}
and let $(x,y,z)$ be the usual
coordinates of $\r^3$. Without loss of generality, we assume that
$\Pi$ is the $xy$-plane, that is, $\Pi=\{z=0\}$ and that $C$ is a
circle of radius $r>0$ centered at the origin. Let
$\vec{a}=(0,0,1)$.

Consider the unit normal vector $N$ given by (\ref{lopez:gauss}). Then
$\langle N,\vec{a}\rangle>0$ on $S$. We will use the notation of
Section \ref{lopez:effect}. First, let  $\alpha$ be the
parametrization of $C$ such that $\alpha'\times\nu=N$. We know
that $\alpha''=-\alpha/r^2$ (independent of the orientation of
$C$) and since $\langle N, \vec{a}\rangle
>0$ along $\partial S$, we have $\alpha\times\alpha'=r\vec{a}$. Then
\begin{equation}\label{lopez:alfa}
\langle\nu,\vec{a}\rangle=\langle N\times\alpha',\vec{a}\rangle=
\langle N,\alpha'\times\vec{a}\rangle=\frac{1}{r}\langle
N,\alpha\rangle.
\end{equation}
The integral equation (\ref{lopez:flux})  gives
\begin{equation}\label{lopez:flux2}
-2\pi r^2 H=\int_{C}\langle\nu,\vec{a}\rangle\ {\rm d}s.
\end{equation}
and thus, equation (\ref{lopez:ba}) is

\begin{equation}\label{lopez:pi}
\int_S\langle N,\vec{a}\rangle\ {\rm d}S=\pi r^2.
\end{equation}

We will need the following result:

\begin{lemma} The function $\langle N,\vec{a}\rangle$ satisfies
\begin{equation}\label{lopez:eq1}
\Delta \langle N,\vec{a}\rangle +|\sigma|^2 \langle
N,\vec{a}\rangle=0,
\end{equation}
where $\Delta$ is the Laplace-Beltrami operator on $S$ and $\sigma$
is the second fundamental form.
\end{lemma}

\begin{proof}
Formula (\ref{lopez:eq1}) holds for any CMC surface. Let
$\textbf{x}:S\rightarrow\r^3$ be an immersion of a surface in
$\r^3$. For any vector field $Y$ of the ambient space $\r^3$, we
consider the decomposition $Y=V+u N$, where $V$ is a tangent vector
field to $\textbf{x}$ and $u=\langle N,Y\rangle$ is the normal
component of $Y$. We consider a smooth variation $(\textbf{x}_t)$ of
$\textbf{x}$ ($\textbf{x}_0=\textbf{x}$) whose variation vector
field is $uN$, that is, $\partial_t(\textbf{x}_t)_{t=0}=uN$. Then
the variation of the mean curvature $H_t$ of the $(\textbf{x}_t)$
changes according to
$$\partial_t(H_t)_{t=0}=\frac12(\Delta u+|\sigma|^2 u)+\langle \nabla H,V\rangle.$$
The first summand in the above equation is the linearization of the
mean curvature operator. See \cite{re1}.

Assume now that $\textbf{x}$ is a CMC surface. Then $\nabla H=0$.
Let us take the vector field $Y=\vec{a}$ whose associated
one-parameter subgroup  generates translations. Thus the mean
curvature  is fixed pointwise throughout the variation, and so,
$\partial_t(H_t)_{t=0}=0$, proving (\ref{lopez:eq1}).

\end{proof}

We follow with the proof of Theorem \ref{lopez:t2}. By  applying   the
Divergence Theorem to the vector field $\nabla \langle
N,\vec{a}\rangle$ and using equation (\ref{lopez:eq1}),  we get
\begin{equation}\label{lopez:eq2}
\int_S |\sigma|^2\langle N,\vec{a}\rangle \ {\rm d}S=\int_{C}\langle
dN \nu,\vec{a}\rangle\ {\rm d}s.
\end{equation}
We study each side of (\ref{lopez:eq2})  beginning with the
left-hand side. The inequality
$|\sigma|^2\geq 2H^2$ holds on any surface, and equality occurs only at  umbilic
points. By using (\ref{lopez:pi}) and
 that $\langle N,\vec{a}\rangle>0$ on $S$, the left-hand side of (\ref{lopez:eq2}) yields
\begin{equation}\label{lopez:eq3}
\int_S |\sigma|^2\langle N,\vec{a}\rangle\ {\rm d}S\geq 2H^2
\int_S\langle N,\vec{a}\rangle\ {\rm d}S=2\pi r^2 H^2.
\end{equation}
We now turn our attention to the right-hand side of (\ref{lopez:eq2}). First, note
$$dN\nu=-\sigma(\alpha',\nu)\alpha'-\sigma(\nu,\nu)\nu.$$
From (\ref{lopez:alfa}), we have
\begin{eqnarray}\label{lopez:eq5}
\sigma(\nu,\nu)&=&2H-\sigma(\alpha',\alpha')= 2H+\langle dN
\alpha',\alpha'\rangle\nonumber\\
&=&2H-\langle N,\alpha''\rangle=2H+\frac{1}{r^2}\langle
N,\alpha\rangle=2H+\frac{1}{r}\langle\nu,\vec{a}\rangle.
\end{eqnarray}
Because $\langle\alpha',\vec{a}\rangle=0$, and using
(\ref{lopez:flux2}) and (\ref{lopez:eq5}), we have
\begin{eqnarray}
\int_{C}\langle dN \nu,\vec{a}\rangle\ {\rm
d}s&=&-\int_{C}\sigma(\nu,\nu)\langle\nu,\vec{a}\rangle\ {\rm d}s
=-\int_C\Big(2H+\frac{1}{r}\langle\nu,\vec{a}\rangle\Big)\langle\nu,\vec{a}\rangle\ {\rm d}s\nonumber\\
&=&4\pi r^2 H^2-\frac{1}{r}\int_C\langle\nu,\vec{a}\rangle^2\ {\rm
d}s.\label{lopez:eq6}
\end{eqnarray}
We use (\ref{lopez:flux2}) and the Cauchy-Schwarz inequality as
follows:
\begin{equation}\label{lopez:eq7}
\int_C\langle\nu,\vec{a}\rangle^2\ {\rm d}s\geq\frac{1}{2\pi
r}\Bigg( \int_C \langle\nu,\vec{a}\rangle\ {\rm d}s\Bigg)^2=2\pi r^3
H^2.
\end{equation}
Then equation (\ref{lopez:eq6}) and inequality (\ref{lopez:eq7}) imply
\begin{equation}\label{lopez:eq8}
\int_{C}\langle dN \nu,\vec{a}\rangle\ {\rm d}s\leq 2\pi r^2 H^2.
\end{equation}
Finally, by combining (\ref{lopez:eq2}), (\ref{lopez:eq3}) and
(\ref{lopez:eq8}), we obtain
$$2\pi r^2H^2\leq \int_S |\sigma|^2\langle N,\vec{a}\rangle\ {\rm d}S=\int_{C}\langle dN
\nu,\vec{a}\rangle\ {\rm d}s\leq 2\pi r^2 H^2.$$
Therefore, we have
equalities in all the above inequalities. In particular,
$|\sigma|^2=2H^2$ on $S$. This means that $S$ is a totally  umbilic
surface of $\r^3$ and so, it is an open of a plane or a sphere.
Because, the boundary of $S$ is a circle, $S$ is a planar disk or it
a spherical cap. In the latter case, $S$ must be the small spherical cap since $S$ is a graph.
This concludes the proof
of  the  theorem.

We point out that the hypothesis that  $S$ is a graph has been only used  in
inequality (\ref{lopez:eq3}) to assert $|\sigma|^2\langle
N,\vec{a}\rangle\geq 2H^2\langle N,\vec{a}\rangle$. However,  the
rest of the proof is valid for any  CMC compact surface bounded by a
round circle. The reader may then try to use the ideas that underlie
the proof of Theorem \ref{lopez:t2} to derive other results. As an example,
we show the following theorem where we replace the hypothesis  that the surface
is a graph by the hypothesis of  non negativity of the Gauss curvature.
Again, we follow the spirit of our work and we do not invoke the Maximum
Principle.

\begin{theorem}
 Let  $S$ be a  compact CMC surface bounded by a round circle $C$. If the Gauss
curvature $K$ is non-negative,
then $S$ is a planar disk or a spherical cap.
\end{theorem}

\begin{proof} As $\langle N,\vec{a}\rangle\leq 1$ holds, we have $K\langle
N,\vec{a}\rangle\leq K$ independent on the sign of $\langle
N,\vec{a}\rangle$. Thus,  $K\langle N,\vec{a}\rangle\leq H^2$. From
(\ref{lopez:pi}) and (\ref{lopez:eq3}) and since
$|\sigma|^2=4H^2-2K$, we have
 \begin{eqnarray*}
\int_S |\sigma|^2\langle N,\vec{a}\rangle\ {\rm d}S&=&
4H^2\int_S\langle
N,\vec{a}\rangle\  {\rm d}S -2\int_S K\langle N,\vec{a}\rangle\ {\rm d}S \\
&\geq& 2H^2 \int_S\langle N,\vec{a}\rangle\ {\rm d}S=2\pi r^2 H^2.
\end{eqnarray*}
The proof then  follows  the same steps as in  the proof of Theorem \ref{lopez:t2}
and we conclude that $S$ is umbilic, and  hence a planar disk or
a spherical cap.
\end{proof}

 We end with a comment. It would be
interesting to have a proof of Theorem \ref{lopez:embedded}, that
is, the boundary version of the Alexandrov theorem, without invoking
the Maximum Principle for elliptic equations, as was done by
Reilly in \cite{re} for the closed case. We do not know a way of
applying our arguments to this question.

\label{lopez:end-art}
\end{document}